\documentclass{amsart}
\usepackage{graphics,verbatim,color,amsthm}

\usepackage[all]{xy} 

\theoremstyle{plain}
\newtheorem{theorem}{Theorem}
\newtheorem{proposition}{Proposition}
\newtheorem{exer}{Exercise}

\newtheorem{definition}{Definition}

\theoremstyle{definition}

\def\R{\mathbb{R}}

\def\C{\mathbb{C}}

\def\bP{\mathbb{P}}
\def\E{\mathbb{E}}

 \newcommand{\dd}[2]
{
{{\partial #1}   \over {\partial #2}}
}

\begin{document}

\title{Blow-up for    realizing homotopy classes in the three-body problem}
 \maketitle
 
 \begin{abstract}
This expository  note describes   McGehee blow-up \cite{McGehee}
in its role as one of the main tools in my recent proof with Rick Moeckel \cite{RM2}
that every free homotopy class for the planar three-body problem can be realized by
a periodic solution.    The main  novelty   is
my use of  energy-balance  
to   motivate the    transformation of McGehee.  Another novelty is 
an explicit description of the blown-up reduced phase space for the planar N-body problem, $N \ge 3$
as a complex vector bundle over  $[0, \infty) \times \C \bP ^{N-2}$ where $r \in [0, \infty)$ measures the size of
a labelled planar N-gon and $[s] \in  \C \bP ^{N-2}$ describes its shape. 
 \end{abstract}
 
\section{Introduction}  
Deleting  collisions endows  the configuration
space of the planar  three-body problem with non-trivial  topology.    Modulo rotations, this   space is homotopic to
a two-sphere minus three points and so has a large set of   free homotopy classes of   loops   
We call these classes the ``relative free homotopy classes.  
 We say that a solution is ``relatively periodic'' if it is periodic modulo rotation, or 
equivalently, if it is periodic in some  rotating frame.
\begin{theorem} [$(RM)^2$ \cite{RM2}] For equal or near-equal masses,  and  angular momenta
$J$ sufficiently small  but nonzero, every relative  free homotopy class for the planar three-body problem is 
 realized by
a relatively periodic orbit for the Newtonian planar three-body problem having energy $-1$ and angular momentum $J$.
\end{theorem}

{\sc Remark.} ``Relative free homotopy classes'' are encoded by ``reduced syzygy sequences'' : periodic lists of 1's 2's and 3's such
as 123232...  where the symbol i indicates  that the three bodies have become  instantaneously collinear with mass i in between j and k
with i,j,k a permutation of the symbols 1,2,3.  See for example (\cite{RM2}) for details.

{\sc History of the theorem.}  A well-known theorem in Riemannian geometry
asserts that on  a compact Riemannian manifold every free homotopy class of loops  is
realized by a  periodic geodesic.   
 Inspired by this   basic  geometric 
fact Wu-yi-Hsiang, in 1996, posed  the question: ``Is every free homotopy
class realized by a (relatively) periodic solution to      the planar Newtonian   three-body equation?''
  For 17  years I tried to use variational methods to  prove that the answer is ``yes'' when the total angular momentum $J$ is zero.
 Finally, at the urging of Carles Sim\'o,  in October of 2014, I gave up on   variational methods and tested the waters of dynamical methods.  
 Almost as soon as I gave up
 I realized that Rick Moeckel had come within epsilon of proving   theorem 1  in the 1980s [ (\cite{M_chaotic}, \cite{M_InfinitelyClose}, and  \cite{M_symbolic})].

My purpose in    this note  is to give an exposition of one  of the principal tools
  in our proof,    McGehee blow-up \cite{McGehee}, and   a sense of how we use this tool to prove the theorem.   
 The main  novelty   is
the use of  energy-balance  
to   motivate the  mysterious transformation of McGehee.  Another novelty is 
an explicit description of the blown-up reduced phase space for the planar N-body problem, $N \ge 3$.   
For further reading on McGehee blow-up we   recommend   Moeckel   \cite{M_near}, pages 222-225  and Chenciner   \cite{Ch_infini}.

\section{Background: Equations and Solutions.}

\subsection{ The Equations.  }

The classical three-body problem demands that we  solve the 
system of non-linear ODEs: 
 \begin{equation}
 \label{N}
\begin{aligned}
m_1 \ddot q_1 & =&  F_{21} + F_{31} \\
m_2 \ddot q_2 & =&  F_{12} + F_{32} \\
m_3 \ddot q_3 & =&  F_{23} + F_{13} .
\end{aligned}
\end{equation}
where
\begin{equation}
\label{force}
F_{ab} = G m_a m_b\frac{ q_a- q_b}{r_{ab}^3} \end{equation}
is the force exerted by mass $m_a$ on mass $m_b$
and
$$ r_{ab} = |q_a -q_b|, q_a \in \R^d,   \text{ and } m_a, G > 0.  $$  
Here $a, b = 1,2, 3$ label the bodies.  The dimension $d$ for us will be $2$. 
(The standard value  is $d=3$. )  
The $m_a$ represent the values of point masses
whose instantaneous positions are $q_a (t)$.  The double dots indicate two time derivatives: $\ddot q = \frac{d ^2 q} {d t^2}$.
 The constant $G$ is Newton's gravitational
constant and is physically needed to make dimensions 
match up.   Being mathematicians, we  can and do set  $G = 1$.

 \subsection{The Solutions of Euler and Lagrange}
\begin{figure}[h]
\label{fig:shape}
\center
\scalebox{0.3}{\includegraphics{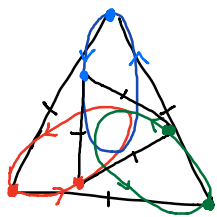}}
\caption{A Lagrange Solution.}  \label{LagrangeSoln}
\end{figure}

\begin{figure}[h]
\label{fig:shape}
\center
\scalebox{0.3}{\includegraphics{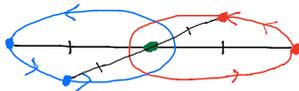}}
\caption{An Euler Solution.}  \label{EulerSoln}
\end{figure}

 The only   solutions to the three-body problem for which we have explicit formulae were  found by Euler \cite{Euler} and Lagrange \cite{Lagrange} in
 the last half of the 18th century.    See figures \ref{LagrangeSoln}, \ref{EulerSoln}. 
 Their solutions are central to our  story. 
 
 For Lagrange's solution, place  the three masses at the vertices of an equilateral triangle and drop them:
 let them go from rest.  They shrink   homothetically shrink towards their common center of mass,
 remaining equilateral at each instant. The solution ends  in finite time in triple collision. This motion   forms half of  Lagrange's triple collision solution.
 To obtain the  other half of Lagrange's solution use time-reversal invariance to continue this solution backwards in time.  In the full solution  the three masses explode out of   triple collision,    reach a   maximum
 size at  the   instant at which we dropped the three masses, and then shrink back to triple collision, staying equilateral throughout.   A surprise  is that the   Lagrange solution 
 works regardless of the mass ratios $m_1:m_2:m_3$.  
  
  For Euler's solutions,   place the   masses on the line in a certain order: $q_i < q_j < q_k$ so as to form a special  ratio $q_k-q_j : q_j -q_i$.  (This special  ratio depends on the mass ratios
 and also  the choice of mass  $m_j$ on the middle and is the root of a fifth degree polynomial whose coefficients depend on the masses.)
   Again   drop them.  They stay on the line as they evolve and  again the similarity class of the (degenerate) triangle  stays constant: 
 this  ratio of side lengths stays constant. 
( In case  the two masses at the ends  are equal  then the special ratio is   $1:1$:  place $m_j$ at  the midpoint
of $m_i$ and $m_k$. )  
 
The solutions just described are part of a   family of explicit solutions discovered by Euler and Lagrange.  
For every one of the      solutions in these families  the similarity class formed by the three masses stays constant in time during the evolution.  Each mass moves on its own
Keplerian conic with the center of mass of the triple as focus, the solutions described i  above being the special case  of
  degenerate (colinear) ellipses.   We   derive these   families   analytically  in section \ref{ELfamily} below.  

All together these solutions form   five   familes.  The corresponding shapes are called ``central configurations''.  The Lagrange solutions count as two, one 
shape for each orientation   of 
a labelled equilateral triangle.   The Euler solutions count as three, one    for each choice of   mass   in the middle.

For almost all (Newtonian) time the solutions  of theorem 1 
are very close to one of the three Euler solutions.  The  Lagrange solutions act as bridges between
various Eulers.

 \eject 
\section{Shape sphere.  Blow-up and reduction, first pass.}
\label{3}

A basic aid to  understanding   the planar three-body problem is the {\it shape sphere},   a two-sphere whose points represent oriented similarity classes of triangles.
At each instant of time  three bodies form the     vertices of a  triangle.  Call two triangles     ``oriented similar'' if one can be brought to the other by a composition of translations, rotations, and scalings. 
The resulting space of equivalence classes forms  the  shape sphere.  See  figure \ref{fig:shape}.
This  sphere has 8 marked points,  the  5 central configurations just described $L_+, L_-,  E_1, E_2, E_3$  and   
 3    binary collision points   labelled $B_{12}, B_{23},  B_{31}$ .  The sphere's equator  represents the   
space of  collinear triangles. The   3 binary collision points,  and    3 Euler central configurations   
lie on this equator, interleaved so as to be alternating.   
\begin{figure}[h]
\label{fig:shape}
\center
\scalebox{0.6}{\includegraphics{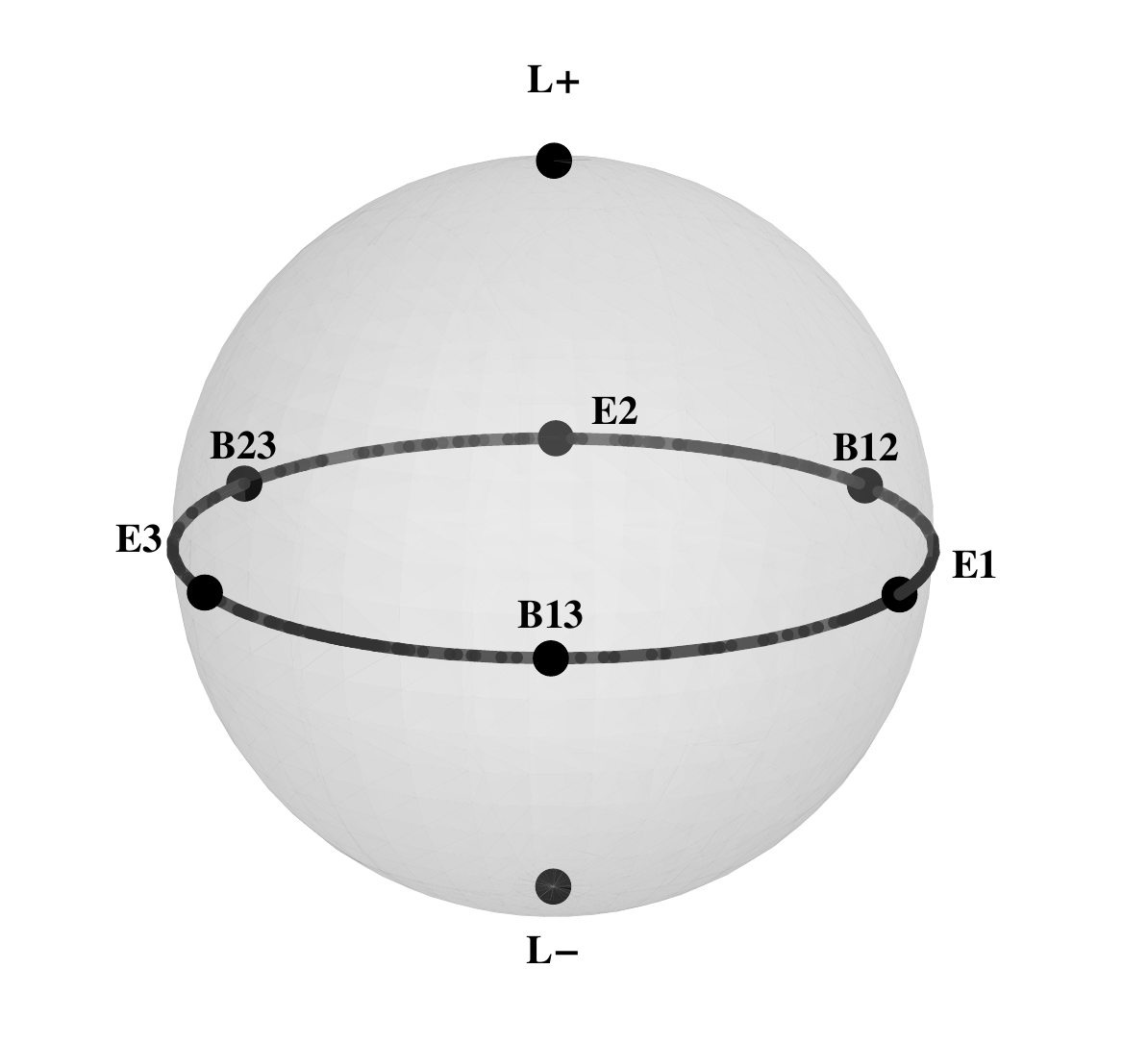}}
\caption{The shape sphere. Lagrange points, Euler points, and collision points marked. The equator consists of collinear triangles. \small{figure courtesy of Rick Moeckel}}  \label{shape sphere}
\end{figure}

The earliest  occuring picture of the shape sphere in the context of celestial mechanics with which I am familiar is  \cite{M_near} .  
You will find  a detailed      exposition of the shape sphere and its relation to   the three-body problem in \cite{me_monthly}. 

We now  summarize how the shape sphere arises out of  the three-body problem.
The configuration space for the 3 body problem, with collisions allowed, is $\C^3$
with  $q = (q_1, q_2, q_3) \in \C^3$ representing the 3 vertices of the triangle - the positions of the 3 bodies.
We have identified $\C$ with $\R^2$ in the standard way: $x + iy \in \C$ corresponds $(x,y)\in \R^2$.
A standard trick  from Freshman physics allows us to restrict the problem to the center-of-mass zero subspace:
$$\E_{cm} = \{ q \in \C^2:  m_1 q_1 + m_2 q_2 + m_3 q_3 = 0 \} \cong \C^2 \subset \C^3.$$
(See  the beginning of section \ref{Reduction} below.) 
In $\E_{cm}$ the binary collision locus become three complex lines which intersect
at the origin $0$. The origin    represents triple collision.  The masses endow $\C^3$ with a canonical metric
called the ``mass-metric'' (eq. (\ref{massmetric})) and relative to that metric the distance from triple colliison is given
by $r$ where 
$$r^2 = m_1 |q_1|^2 + m_2 |q_2|^2 + m_3 |q_3|^2 .$$
(See  eq. (\ref{size}).) 
Take the sphere 
$$\{ r = 1 \} := S^3 \subset \C^2 \cong \E_{cm}.$$
Because the  three-body equations are invariant under rotations they descend to ODEs on the  quotient
of  $\C^2 = \E_{cm}$ by the group $S^1$ of rotations.  This quotient space is
topologically an $\R^3$.  To understand this quotient   note that the rotation   action leaves $r$ unchanged 
but moves points on $S^3$ around according to $(Z_1, Z_2) \mapsto (uZ_1, uZ_2)$, $u \in S^1 \subset \C$.
(Here $Z_1 , Z_2$ are any complex linear coordinates for $\E_{cm}$.)
This is the circle action used to form the Hopf fibration: 
$$Hopf: S^3 \to S^3 / S^1 = S^2 = \text{ shape sphere }.$$
Points of the quotient $\R^3$  represent oriented congruence classes of triangles:  planar triangles modulo translation and rotation,
 but not scaling.   Express  $\R^3$  in   spherical  coordinates $(r, s),  s \in S^2 $.  Then  the origin
 $r=0$ corresponds to   triple collision.   A point $s$ on the sphere represents a  ray $rs,  r\ge 0$ of triangles all having the same shape.  
 The collision locus ${\mathcal C} = \{ r_{12} = 0 \text{ or } r_{23} = 0 \text{ or } r_{31} = 0 \}$  is then represented by the three rays corresponding to the three binary collision points $B_{12}, B_{23}, B_{31} \in S^2$.

Newton's equations  break down at triple collision  $r=0$.  McGehee
blow-up is a  change of variables  ( equations (\ref{McGehee})) 
which converts  Newton's equations to a system of ODEs which is well-defined when   $r=0$.  The locus 
$r=0$ in the new variables  is called ``the collision manifold'' and forms a bundle over the shape sphere.
The blown-up system of ODEs has exactly 10 fixed points, all   on the   collision manifold ,
with a   pair of fixed points lying over each of the five central configurations.   For a chosen central configuration,  one   element of the     pair corresponds to the  homothetic arc  incoming to triple collision,
as in our original description of the  Lagrange solution,  while the  other element of the pair corresponds  to the initial segment of that solution which explodes out from  triple collision.

The 10 fixed points  on the collision manifold    have stable and unstable manifolds,  parts of which stick out of the collision manifold, and which 
 intersect in complicated ways, as per the Smale Horseshoe and heteroclinic tangles.   See figure \ref{web}.
Moeckel   investigated these manifolds and their relations  in    seminal works \cite{Moeckel1}, \cite{M_InfinitelyClose},
 \cite{M_heteroclinic}, \cite{M_near}, \cite{M_chaotic},   and \cite{M_symbolic}
where he  proved existence of ``topological heteroclinic tangles''  between them. 

One finds the following abstract   graph 
\vskip .9cm

\begin{figure}[h]
\scalebox{0.4}{\includegraphics{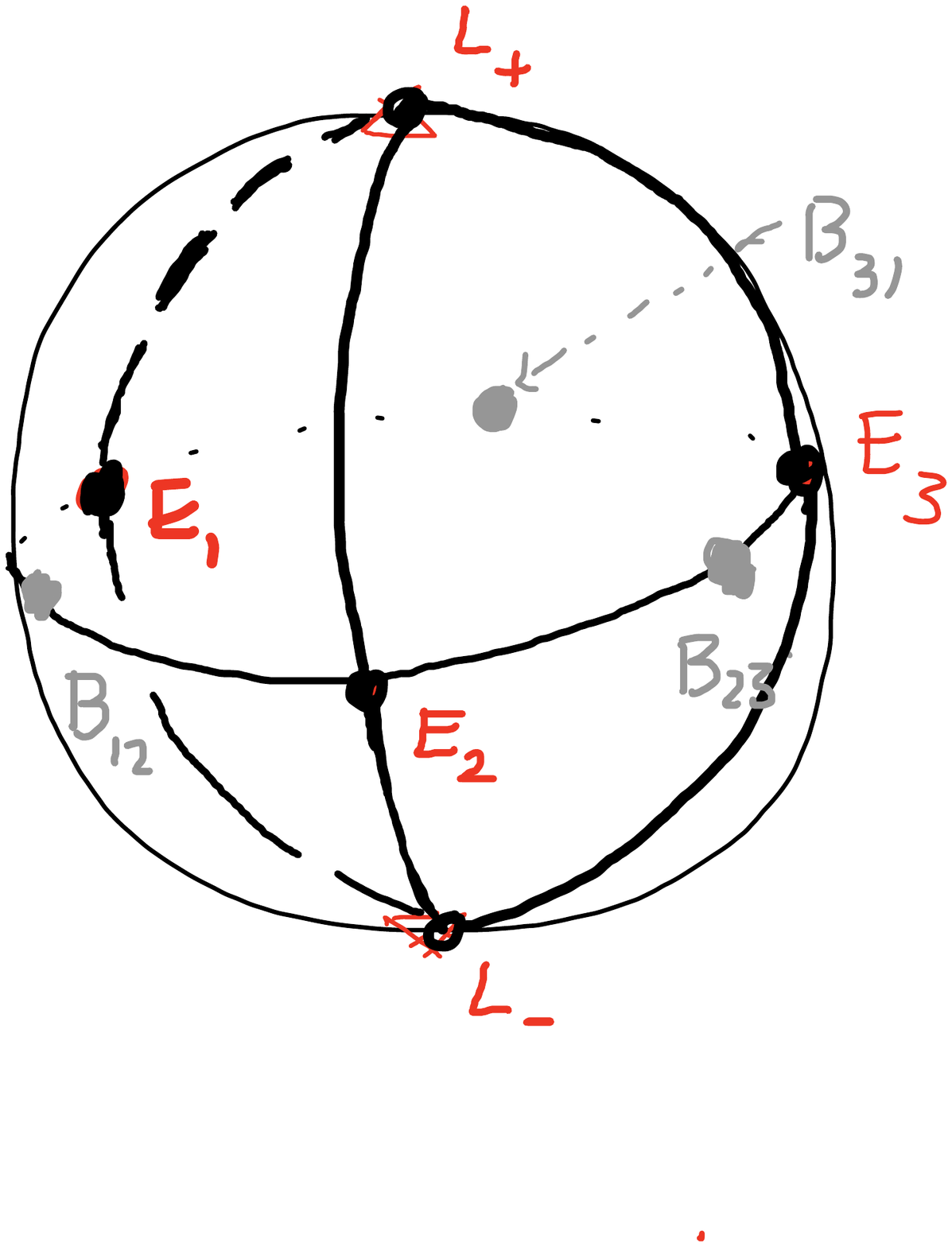}}
\caption{The concrete graph, embedded in the shape sphere.}  
\label{graphWSphere}
\end{figure}

\begin{figure}[h]
\scalebox{0.4}{\includegraphics{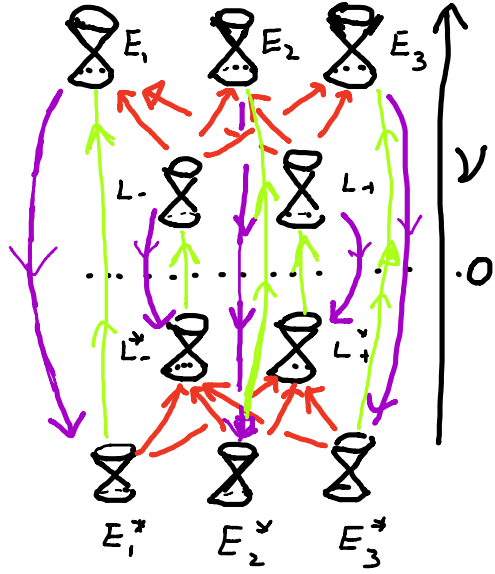}}
\caption{The equilibria arising upon blow-up and relations between their stable and unstable manifolds. The purple and green arrows
are the rest cycles described in figure \ref{fig:ccplots}.  }  
\label{web}
\end{figure}

\hskip 1.5cm 
$$ $$
\begin{picture}(0,0)

\put(30,0){\circle{2}}

\put(32,-4){$L_-$}

\put(30,0){\line(0,1){25}}

\put(30,0){\line(-1,1){26}}

\put(30,0){\line(1,1){26}}

\put(30,25){\circle{2}}

\put(32,26){$E_2$}

\put(5,25){\circle{2}}

\put(10,26){$E_1$}

\put(55,25){\circle{2}}

\put(57,26){$E_3$}

\put(5,25){\line(1,1){26}}

\put(30,25){\line(0,1){26}}

\put(55,25){\line(-1,1){26}}

\put(30,50){\circle{2}}

\put(32,52){$L_+$}

\label{pic} 
\end{picture}
$$ $$
 in several of these  papers of Moeckel (\cite{M_chaotic}, p. 53, Theorem $1^{\prime}$, and   \cite{M_symbolic} whose  Figure 2 becomes our graph
 after  deleting   the vertices labelled  with Bs  (for    binary collision) as well   edges incident to them).   
Moeckel's theorem in  \cite{M_chaotic}, based on the intersections between   stable and unstable manifolds of the 10 fixed points, asserts that all paths in this graph are ``realized'' by solutions to the three-body problem provided the angular momentum, energy
and masses are as per theorem 1.  
Embed this graph in the  shape sphere as indicated by figure \ref{graphWSphere}.  Call the   embedded graph
the ``concrete connection graph''.   

 The dynamical relevance  of the concrete connection  graph has to do with the Isosceles three-body problem. 
 When two  of the masses are equal, say $m_1 = m_2$, then the isosceles triangles  $r_{13} = r_{23}$
form an  invariant submanifold
of the three-body problem whose dynamics is called the ``Isosceles three-body problem''.  These Isosceles triangles
form a great circle passing   through both  Lagrange points,   the binary point $B_{23}$ ,  and the Euler point $E_1$.   If all three masses are equal  we have
   have  three Isosceles subproblems represented by three great circles on the shape sphere. Take one-half of each  great circle,  namely that half whose endpoints are the two Lagrange points
  and which contains   the Euler point. In this way we form the concrete connection graph in which the  edges are parts of the Isosceles great circles.    
  
  Observe that the shape sphere minus the three binary collision points retracts onto the concrete connection  graph. 
Theorem 1   follows immediately from this observation and Moeckel's theorem refered to above,   once we know that the realizing   solutions of Moeckel's theorem,
 projected onto the shape sphere,  stay $C^0$-close   to   corresponding edges in the  concrete  connection graph.   For a few more details see the final section of
 this article.

\vskip .2cm 
\section{Set-up.  Blow-up}

It is no more work
to perform the blow-up for the N body problem in d-dimensional Euclidean space, rather than   our special case
of the three-body problem in the plane.
The d-dimensional N-body equations are:    
\begin{equation}
m_a \ddot q_a  = \Sigma_{b \ne a}  F_{ba} \qquad, q_a   \in \R^d
\label{eq:N}
\end{equation}
with the forces $F_{ba}$ as  above.

\subsection{Metric Reformulation.}

Let 
$$\E= (\R^d) ^N$$
denote the N-body configuration space.  Write points of $\E$
as  $q = (q_1, \ldots , q_N)$ and think of the points  as the N-gons in d-space.    The masses endow $\E$ with an   inner product,
called the   {\it mass inner product} :  
\begin{equation}
\langle q, v
\rangle = \Sigma m_a q_a \cdot v_a
\label{massmetric}
\end{equation}
so that the standard     kinetic energy is given by  
\begin{equation}
\label{kinetic}
K = \frac{1}{2}\langle \dot q, \dot q \rangle.
\end{equation}
 Let $\nabla$ be the   gradient associated
to this metric:
$df_q (v) = \langle \nabla f(q), v \rangle$, so that
$(\nabla f)_a = \frac{1}{m_a} \dd{f}{q_a}$.
Then the N-body equations take the simple form 
\begin{equation}
\ddot q = \nabla U (q)
\label{eq:Nd}
\end{equation}
where $U$ is the {\it negative} of the standard potential $V$:  
 \begin{equation}
\label{potential}
 U =- V =  \Sigma_{a < b}  \frac{m_a m_b }{r_{ab}},
 \end{equation}
 the sum being over all distinct pairs
$a, b$.   As is well known, the   total energy is 
conserved (i.e constant) along solutions 
$$H = K - U = K + V .$$

We use 
\begin{equation}
\label{size}
r = \sqrt{\langle q, q \rangle}
\end{equation}
to  measure the size
of our configuration $q = (q_1, \ldots , q_N)$.
Lagrange proved that
$$r^2 = \Sigma_{a < b}  m_a m_b r_{ab} ^2 / \Sigma m_a$$
provided we are in center of mass coordinates: $\Sigma m_a q_a = 0$.
Then
$$r = 0 \Leftrightarrow \hbox{total collision: all masses coincide}$$
while 
$$U   = \infty  \Leftrightarrow \hbox{some collision: some pair of masses coincide} .$$
\begin{exer} Use the metric reformulation of Newton's equations  eq (\ref{eq:N}), the fact that $U$ is homogeneous of degree $-1$
and   Euler's identity for homogeneous functions to  derive the ``virial identity'', also known as the Lagrange-Jacobi identity:
$d ^2 (r ^2) /dt ^2 = 4H + 2 U$. Also show that $2H + U = H + K$.
\end{exer}

\subsection{ McGehee transformation via Energy Balance.}

The key property of  the potential energy 
$-U$,   as far as McGehee's transformation is concerned, is that it is homogeneous of degree $-1$: $U (\lambda q) = \lambda ^{-1} U(q)$,
or $$q \mapsto \lambda q \implies U \mapsto \lambda^{-1} U.$$
Our guiding principle in deriving the McGehee blown-up equations is to require that the kinetic energy $K$ scale the same   as   potential energy.
  so that the total energy balance under scaling and thus  has a scaling law.  We call this principle ``energy balance''. Then $K \mapsto \lambda^{-1} K$ and since
 $K$ is quadratic in velocities $v$   velocities  must scale by
 $$v = \dot q \mapsto \lambda^{-1/2} v.$$

How  must time scale?  
Since $dq \mapsto \lambda dq$  
and  $v = dq/ dt$, we see that for a power law scaling  $dt \mapsto \lambda ^{a} dt$ to yield 
  $v \mapsto \lambda ^{-1/2} v$ we must have $a =3/2$.  Summarizing, our space-time  scaling law must be 
\begin{equation}
q \mapsto \lambda q,  dt \mapsto \lambda^{3/2}dt, 
\label{spacetime}
\end{equation} 
which   induces the desired scalings
$$ v \mapsto \lambda ^{-1/2} v ;   (U, K , H) \mapsto \lambda ^{-1} (U, K, H)$$


 \begin{exer}   
  Show that $q(t)$ solves (\ref{eq:Nd}) if and only also $q_{\lambda} (t) :=  \lambda q (\lambda^{-3/2} t)$ solves (\ref{eq:Nd}).  Explain how the exponent  $-3/2$
   in this   transformation-of-paths formula arises from the   $+3/2$ in the time part of the scaling law of  eq   (\ref{spacetime}). 
   \end{exer}

 McGehee's genius    was to rewrite Newton's equations, as  much as is possible, in   scale invariant terms.
We cannot completely get rid of scale, but we can encode  scale in   the single size  variable  $r = \sqrt{\langle q, q \rangle} $ introduced earlier and through which  we   remove scale from the remaining  variables: 
\begin{eqnarray}
q & = & rs \\
v & = & r^{-1/2} y \\
dt & = &  r ^{3/2} d \tau 
\label{McGehee}
\end{eqnarray}
 These relations define the
 {\it McGehee transformation} $(q,v; t) \mapsto (r, s, y; \tau)$.
Observe that $s$ lies on the unit sphere $r =1$ in the configuration space,
$$s \in S = S^{dN-1} = \{r =1 \}   \subset \E$$
 so that $(r,s)$ are spherical coordinates on $\E$.  We sometimes refer to    $s$ as  the shape 
 of the configuration $q$.  
 \begin{exer} 
 Write $^{\prime}$ for $\frac{d}{d \tau} = r^{3/2} \frac{d}{dt}$.
Show that   McGehee's transformation  transforms     Newton's equations (\ref{eq:Nd})
to the equations
\begin{equation}
\begin{aligned}
r^{\prime}& =  r \nu \\
s^{\prime} & =  y - \nu s  \\
y^{\prime} & =   \nabla  U (s) + \frac{1}{2} \nu y
\end{aligned}
\label{blowup}  
\end{equation}
 where $\nu = \langle s ,y \rangle$.
 These equations are the McGehee blown-up equations.
\end{exer}

In the last equation $\nabla U (s) \in \E$ is the same gradient as in Newton's equation (\ref{eq:Nd}),
only  that restricted to  points $s$ of the sphere $\{r =1\}$.
The blown-up  equations are analytic and extend analytically to the total collision manifold $r =0$.  For $N>2$ the equations still have singularities
due to partial collisions eg $r_{12} =0$, at which $\nabla U(s)$ still blows up.  
\begin{definition}  The ``extended collision manifold''  is the locus $r = 0$ 
for the blown-up phase space $[0, \infty) \times S \times \R^{dN}$
of  McGehee. 
\end{definition}  


The first of the three blown-up ODEs  asserts that the 
extended collision manifold is  an invariant submanifold.    On the extended collision manifold the flow is non-trivial , as a glance at
the last two equations shows.     Away from the extended collision manifold, the blown-up equations
are equivalent to Newton's equations.   
What have we gained by adding this collision manifold?

\subsection{Equilibria!  }   

The first thing one learns  in a class in dynamical systems is to look for equilibria.  
But  Newton's equations have no equilibria! N stars cannot just sit there, still, in space. 
Adding the extended  collision manifold thru blow-up  introduces  equilibria. 
When $N =3$ these equilibria
correspond   to the  solutions of Euler and Lagrange described above.  In the
general case the equilibria  correspond to ``central configurations''.  (Proposition \ref{ccs} below.)  

{\it Finding the equilibria.}  From the first of the blow-up equations (\ref{blowup}) 
we see that at an equilibrium  must lie on the extended collision manifold  $r = 0$ (consistent with
wha we just said about ``stars cannot just sit there'').
Plugging the second equilibrium equation $0 = y- \nu s$ into the third equation of the blown-up equations (\ref{blowup})  yields  the ``shape equation'' 
$$\nabla U (s) = - \frac{1}{2} \nu^2 s.$$
Taking the inner product of both sides of the shape equation with $s$ and using Euler's identity for homogeneous functions
yields $$U(s) = \frac{1}{2} \nu^2$$ or $\nu = \pm \sqrt{ 2 U(s)}$.
Now      the   gradient of the function $r^2$   at a point  $s \in S$
is $2s$ so  that we can rewrite the  shape equation  as  
 $$\nabla U (s) =  c \nabla (r^2) \qquad,  c= - \frac{1}{4} \nu^2 = -\frac{1}{2} U(s). $$
 Think of $c$ as a Lagrange multiplier.   We   have proved that the  
 shape $s$ of an equilibrium configuration  must be  a ``central configuration'' where:
\begin{definition} A central configuration    is
a shape $s \in S$ which is a critical point of $U$ restricted to the sphere $r =1$.
\end{definition}
Conversely, for  each central configuration shape  $s_{cc} \in S$ we 
 obtain an equilibrium point $(r,s, y) = (0, s_{cc}, y)$ by setting  $y = \nu s_{cc}$.
with   $\nu = \pm \sqrt{2 U(s_{cc} )}$..
  We have established 
 \begin{proposition} 
 \label{ccs} Equilibria of the blown-up equation are in 2:1 correspondence  with
 central configurations. This correspondence associates to a given central configuration
 $s_{cc}$ the two equilibria $(r,s,y) = (0, s_{cc}, \nu s_{cc})$, 
 with $\nu = \pm\sqrt{2 U(s_{cc})} s_{cc})$.
 \end{proposition}
 
 There is another way to arrive at central configurations in keeping with our original
 discussion of the Euler and Lagrange solutions. Make the
 ansatz:
 \begin{equation}
 q(t) = \lambda (t) s
 \label{ansatz}
 \end{equation}
 where $\lambda (t)$ is a time dependent scalar and $s \in S$ is constant.
 \begin{exer}
 \label{STAR}
 Show that the ansatz (\ref{ansatz}) satisfies Newton's equations if and only if
 $s$ is a central configuration and $\lambda(t)$ satisfies the ``Kepler problem'':
 \begin{equation}
 \ddot \lambda = -\mu  \lambda / |\lambda |^3 \qquad \mu = \frac{1}{2} \nu^2= U(s)
 \label{Kepler}
 \end{equation}
 \end{exer}
 All solutions of the one-dimensional  Kepler problem end in collision:  $\lambda = 0$.
 This value of $\lambda$ in the ansatz corresponds to total collision.  
 The ansatz (\ref{ansatz}) with real scalar $\lambda (t)$ yields the 
 Lagrange and Euler solutions to the three body equations which  we  first described above  by ``dropping'' bodies.

 \subsubsection{The  Euler and Lagrange family.  Planar problems.}
 \label{ELfamily}  Assume we are in the  planar case so that  $d =2$.   Identify $\R^2$ with $\C$: $(x,y) \to  x+iy$
 so that $\E = \C^N$ and so that   complex scalar multiplication of $s = (s_1, \ldots , s_N) \in  S \subset   \C^N$  by $\lambda \in \C$
corresponds to scaling the N-gon $s$  by the factor  $|\lambda|$ while rotating it by $Arg(\lambda)$.
\begin{exer} 
\label{equiv}
Show that $\nabla U(\lambda q) = \frac{\lambda}{|\lambda|^3} \nabla U (q)$
\end{exer}
\begin{exer} 
\label{planarccs}
Use exercise \ref{equiv} to show that exercise   
\ref{STAR} also  holds in the case of the planar N-body problem,  with now    $\lambda(t) \in \C$ a {\it complex scalar}.
\end{exer}
The solutions  of exercise \ref{planarccs}
 are motions in   which all N bodies move  ``homographically'' : meaning by scaling and rotating.
 Since $\lambda (t)$ parameterizes a  conic, each body moves along a homographic conic $\lambda (t) s_a$, $a =1, 2 \ldots, N$.    
 We  now   have, for  each {\it planar} central configuration, a  family of solutions parameterized by  the complex solutions to eq. (\ref{Kepler}), 
 and varying from total collision solutions  when $\lambda (t) \in \R$   to circular motions
 $\lambda (t) = e^{i \omega t} \in S^1 \subset \C$.   For fixed energy $h$ we can think of the parameter of the family as the angular momentum $J$
 discussed below, with $J=0$ being total collision and the maximum or  minimum value of $J$ being the circular motion.

 \subsubsection{An open problem}
 
 The potential $U$ is invariant under rotations and translations.  Consequently, the
 central configurations as we defined them are not isolated, but come in families. 
 
{\bf  Is  the set of  central configurations, modulo rotatations and translations,   a finite set?  }
  This problem is attributed to Chazy \cite{Chazy}.  See Albouy-Cabral \cite{Albouy_Cabral} and
  for perspective and  a recent survey.
 
{\it  What is known. Some History. }  $N=3$:  Euler and Lagrange had established  the complete list of central configurations as described here.
 $N=4$: Albouy \cite{Albouy} classified the central configurations in the case of  4 equal masses two centuries two decades and a few years after Euler and Lagrange .  One of his
 main achievements  was to show that in the equal mass  case the 4-body central configurations
 all have a  reflectional symmetry. 
 Eleven  years after Albouy's work   Hampton and Moeckel  \cite{Hampton} proved that the central configurations are finite (less than $1856$ XXX in number!) 
$N=5$: In 2012  Albouy and Kaloshin \cite{Albouy_Kaloshin}   proved that for   $N=5$ and away from  an algebraic surface in
  the parameter space $\R \bP^4 _+$ of mass ratios, the number of central configurations is finite.  
  In 1999 Roberts \cite{Roberts} constructed examples for $N=5$ with one of the five masses negative in which  the set of central configurations 
  is infinite, underlining the subtlety of the problem.

 \subsection{Linear and angular momentum}
 
 Besides  energy, the only  known constants of the motion
 for the general N-body problem are the components of the  linear momentum
 $$P = \Sigma m_a v_a$$
 and the angular momentum
 \begin{equation}
 J(q,v) = \Sigma m_a q_a \wedge v_a
 \label{angmom}
 \end{equation}
 These momenta are intimately connected to the fact that the group $G$
 of rigid motions acts by symmetries of Newton's equations.
 \begin{exer}
 \label{group_perp}
 $v \in \E$ is orthogonal to the $G$ orbit thru $q \in \E$
 if and only if $P(v) = 0$ and $J(q,v) = 0$
 \end{exer}
 
 \subsection{Center of mass frame} 
 
 A well-know argument using Galilean symmetry and  found in essentially any  
introductory physics text allows us to  
suppose that all our solutions satisfy $P = 0$
and 
$$\Sigma m_a q_a = 0$$
In this case we say that we are in ``center of mass frame''
and 
we set
$$\E_{cm} = \{ q \in \E:  \Sigma m_a q_a = 0 \} \cong \R^{d(N-1)}.$$
The infinitesimal generators of the translation action are 
the vectors $q = (c,c, \ldots , c) = ``c \vec 1'',   c \in \R^d$. Now
 $\E_{cm}$ is precisely the orthogonal complement to the subspace of  vectors of
 the form $c \vec 1$.
 This space of   vectors corresponds to the generators of  the translation group,  
 or alternatively to the space  of all   total collision configurations.  It follows
that $\E_{cm}$ realizes  the quotient of $\E$ by translations
and that in $\E_{cm}$   only one point    represents  total collision:  the origin. 

We can  go to center of mass frame before or after blow-up, the result is the same, namely
the system of ODEs (\ref{blowup})  (with poles on the partial collision locus)  but restricted to the the subset of   variables 
\begin{equation}
(r, s, y) = [0, \infty) \times S_{cm} \times \E_{cm}
\label{blowupphasespace}
\end{equation}
where
$$S_{cm} = \{q \in \E_{cm}:  \langle s, s \rangle  = 1 \}\cong S^{d(N-1) -1}.$$
  
 \subsection{Energy-momentum level sets and the Standard   Collision Manifold.} 
 
Because energy and angular momentum are invariant as we flow according
to Newton, by  fixing  their values $h$ and $J_0$ we obtain invariant submanifolds
of phase space: 
  $$M^{int}(h) = \{H = h , r> 0 \}$$
  and
  $$M^{int} (h, J_0) = \{H = h , r> 0 \}$$
 Energy and angular momentum are not defined at $r=0$ so we have excluded $r =0$.
 Set
\begin{equation}
\label{energymom}M(h) =\text{Closure}(M^{int} (h)), \qquad
 M(h,J_0) = \text{Closure}(M^{int} (h, J_0)),
\end{equation}
the  closure  being within within the blown-up phase space. 
 We will  need to understand  the boundaries of these
spaces, which is their intersection with the extended collision manifold $r =0$;
in other words we must understand how these invariant  submanifolds approach  the extended collision manifold
 $\{r =0 \}$ as $r \to 0$.

The following notation will be useful in this endeavor.
\begin{definition}  {\bf [Notation]}  For $F = F(q,v) $ a homogeneous function on $\E \times \E$ 
write $\tilde F$ for the scale-invariant version of $F$    achieved by multiplying $F$ by $r^{-\alpha}$
where $\alpha$ is the degree of homegeneity of $F$ with respect to   our weighted scaling.
Thus:  $F(q, v)  = r^{\alpha} \tilde F(s, y)$. 
\end{definition}
According to ``energy balance'' both the potential energy, kinetic energy, and total energy are homogeneous of degree $-1$. Thus  
$$ \tilde U (s) = r U(q)$$
where $\tilde U$ is homogeneous of degree $0$ and can be viewed as a function on the sphere $S_{cm}$. And 
\begin{equation}  \tilde H = r H
\label{energy}
\end{equation}
where 
$$\tilde H(s,y) = \frac{1}{2} \langle y, y \rangle - U(s)  = \tilde K(y) - \tilde U(s)$$
and $\tilde K, \tilde U$ are homogeneous of degree $0$.
The angular momentum  is 
homogeneous of degree $1/2$ so
that
\begin{equation}
\label{blownupmomentum}
J = r^{1/2} \tilde J(s,y)
\end{equation}
where $\tilde J$ is scale invariant and equals $\Sigma m_a s_a \wedge y_a$.

If follows immediately from   eq (\ref{energy})  
that  
 $$ \partial (M(h)) = \{ \tilde H = 0,  r =0 \} . $$ 
 while using in addition eq (\ref{blownupmomentum}) we see that 
 $$\partial (M(h,0)) = \{ \tilde H =0, \tilde J = 0, r = 0 \}.$$
 These are basic important submanifolds so we give them separate names.
\begin{definition} 
\label{fullcollision}The full collision manifold is $M_0 = \{ \tilde H = 0,  r =0 \} $.  
\end{definition}
 \begin{definition} 
 \label{stdcollision} 
 The   ``{\it standard collision manifold}'' is the locus
$$C:= \{r= \tilde H = \tilde J = 0 \} .$$ 
\end{definition}
Thus the extended collision manifold contains the full collision manifold $M_0$ which in turn
contains the standard collision manifold $C$.
 The equilibria all lie on $C$. 
  Another    reason  for the importance of the   standard collision manifold $C$
is  a  theorem of Sundman. 
\begin{theorem} (Sundman) If $r \to 0$ along an honest  solution, then $J=0$ for that solution
and hence that solution tends to $C$ as $r \to 0$.  Moreover, the solution tends to the subset of
equilibria within $C$. 
\end{theorem}
Here we are using the hopefully obvious
\begin{definition}  
 An ``honest solution'' to the blown-up equations is 
a solution such that $r >0$.  
\end{definition} 
\noindent The honest solutions   are just the reparameterizations
of solutions to our original  Newton's equations according to the blown-up time.

{\bf Remark.} 
The standard collision manifold $C$
 is the space most authors    refer to when they speak of the ``collision manifold''
for the N-body problem.   Chenciner  (see also \cite{Chenciner1}) argues that the standard collision  manifold is the dilation quotient
of the   N-body phase space. 

\begin{exer} Use eq (\ref{blowup}) to show that 
$$\frac{d}{d \tau}{\tilde H} = \nu \tilde H.$$
$$\frac{d}{d \tau}  \tilde J  = -\frac{1}{2} \nu \tilde J$$
hold everywhere on the blown-up phase space.
\end{exer}
It follows from this exercise that $\{ \tilde H = 0\}$ and $\{\tilde J = 0\}$ are invariant manifolds,
as are $M_0$ and $C$.

\section{Quotient by Rotations.}
\label{Reduction}

 Newton's  equations  and their McGehee blow-ups  (eq \ref{blowup}) are invariant
 under the group  $G$ of rigid motions and so descend to ODEs on the quotient space
 of their phase spaces by  $G$.    Working on this quotient instead of the original
 helps our intuition enormously, especially in the case $N=3$ and $d=2$. 
 We describe the quotient and some aspects of the quotient flow.
  
   The group $G$ of rigid motions is the product of  two subgroups, the translation group and the rotation group. We  have  
 already formed the quotient of phase space by translations when we went to  center-of-mass frame, i.e. by  restricting to $s, y \in \E_{cm}$.
 To form the remaining quotient by rotations it is much cleaner to restrict to the planar case $d=2$.  {\it Henceforth we assume that we
are working with the planar N-body problem, $d =2$.} We identify $\R^2$ with $\C$ as before. 
Thus $\E \cong \C^N$ and $\E_{cm} \cong \C^{N-1}$.
 Represent rotations as unit complex scalars $u \in S^1 \subset \C$ acting
on  $(q, v) \in \E_{cm} \times \E_{cm}$
by $(q, v) \mapsto (uq, uv)$ and on   McGehee coordinates
by $(r, s, y) \mapsto (r, us, uy)$.     
\begin{definition}
The blown up reduced phase space in the planar case is the quotient of the
blown-up center of mass phase space $[0,\infty) \times S_{cm} \times \E_{cm} \cong [0, \infty) \times S^{2N-3} \times \C^{N-1}$
by the group of rotations.  Upon deleting the collision locus ${\mathcal C}$ we denote this quotient by 
$${\mathcal P}_N = ([0, \infty) \times (S^{2N-3}\setminus {\mathcal C})  \times \C^{N-1})/S^1.$$
\end{definition}

Momentarily  
forget the velocities $v$ or $y$, and the deletion of the collision locus ${\mathcal C}$  in trying to understand the quotient. 
 The  circle action sends a  blown-up configuration $(r, s)$ to $(r, us)$,  
$s \in S_{cm}$.  So we   need to understand the quotient of 
the sphere $S_{cm} =S^{2N-3}$ by this action of $S^1$.
It is well known that this quotient   $S_{cm}/ S^1$ is isomorphic to the complex
projective space $\C \bP ^{N-2}:= \bP (\E_{cm})$ with the projection map $S_{cm} \to S_{cm}/S^1$
being the Hopf fibration.  Hence the quotient of the $(r,s)$
by $S^1$ yields $[0, \infty) \times \C \bP^{N-2}$.

To better understand the  meaning of  points of  $\C \bP ^{N-2}$, 
work with  $q \in \E_{cm}$ instead of $q =s \in S_{cm}$, insisting only that $q \ne 0$ and  
 now allowing the scalar  $u$ to vary over the larger group $\C^* \supset S^1$ of   {\it all} nonzero complex numbers.
The resulting quotient is   well-known to be
 $(\C^{N-1} \setminus \{0 \})/ C^* = \C \bP^{N-2}$.   The action of $u \in \C^*$ on $q \in \C^{N-1} \setminus \{0\}$
is precisely the action of rotating and scaling the (centered)  N-gon $q$.   
\begin{definition}  The projective space $\C \bP^{N-2}$ just constructed is  called {\it shape space}.
 Its points represent {\it oriented similarity
classes} of planar N-gons.  
\end{definition}

We have realized the configuration part of the quotient after blow-up as $[0, \infty) \times \C \bP^{N-2}$
where $\C \bP^{N-2}$ is the shape space.   {\it When $N=3$ the shape space is the shape sphere
described above.}

{\bf Collision locus.}  

The condition that a configuration $q = (q_1, \ldots , q_N)$ represent a collision  is
that $q_a = q_b$ for some $a \ne b,  1 \le a, b \le N$. 
This condition is complex linear when viewed in  homogeneous coordinates
$[q_1, q_1, \ldots , q_N]$ and
so defines a complex hyperplane,  a $\C \bP ^{N-3} \subset \C \bP ^{N-2}$.
There are   are $N \choose 2$  pairs $(a,b)$ and so we have to delete $N \choose 2$ hyperplanes from 
our shape space.    The union of these hyperplanes, viewed projectively,  is the collision locus: 
$${\mathcal C} = \{[q] = [q_1, q_2, \ldots , q_N]  \in \C \bP ^{N-2}:  q_a =q_b \text{ some } a \ne b \}.$$
We use the same symbol for the collision locus before or after quotient.

{\bf Accounting for velocities.} 
In the last few paragraphs above we dropped  the velocity  $y$.
The quotient map $(r, s) \mapsto (r, [s])$
from  $[0,\infty) \times S_{cm} \to [0, \infty) \times \C \bP ^{N-2}$
expresses $[0,\infty) \times S_{cm}$ as a principal $S^1$ bundle over
$ [0, \infty) \times \C \bP ^{N-2}$.  

Now include the velocity $y$. The quotient procedure with $y$ included   is
precisely the procedure used to construct an   associated vector bundle   to a principal bundle.
Realizing this, we see that  the quotient ${\mathcal P}_N$ is a complex   vector bundle over $[0,\infty) \times (\C \bP^{N-2} \setminus {\mathcal C} )$
whose  rank is  $N-1$ -- the fiber being  coordinatized by  $y \in \E_{cm}$.
What is this vector bundle?  

\begin{proposition} 
\label{blowupQuotient}  $${\mathcal P}_N = [0,\infty) \times T (\C \bP^{N-2} \setminus {\mathcal C} ) \times \R^2$$ 
as a vector bundle over $[0,\infty) \times (\C \bP^{N-2} \setminus {\mathcal C} )$.    The final
$\R^2$ factor    is coordinatized by $(\nu, \tilde J)$
where $\nu = \langle s, y \rangle$ represents the time rate of change of size and  where $\tilde J = \langle is, y \rangle$ is 
also equal to $r^{-1/2} J$ off of $r = 0$ where $J$ is the usual  total angular momentum
of the system. 
The fiber variable  tangent   to   shape space $\C \bP^{N-2}$   represents  ``shape''   velocity.  
\end{proposition}

In the  case of $N=3$ we have    $\C \bP^{N-2} = \C \bP ^1 = S^2$, the shape sphere previously discussed in section \ref{3}.
Then
$${\mathcal P}_3 = [0, \infty) \times  T(S^2 \setminus {\mathcal C}) \times \R^2 =  [0, \infty) \times \R \times (S^2 \setminus {\mathcal C}) \times \R^3 $$
where 
$${\mathcal  C} = \{ B_{12}, B_{23}, B_{31} \}$$ 
is the set of three binary collision points.
 
 \subsubsection{Velocity (Saari) decomposition.}
 Passing thru a configuration $q \in \E_{cm}$ we have two group-defined curves: the scalings $\lambda q, \lambda \in \R$ of $q$
and the rotations $uq,  u \in S^1$ of $q$.    The tangent spaces to these curves are orthogonal, 
and  together with the  orthogonal complement of their span they 
define a geometric splitting of $T_q \E_{cm} = \E_{cm}$
\begin{eqnarray}
T_q \E_{cm}  & =& \text{(scale)} + \text{ (rotation)} + \text{ (horizontal )}
\\
&=& \R q \hskip .3cm \oplus \hskip .3cm i \R q \hskip .3cm \oplus \hskip .3cm \{v: J(q,v) = 0 , \nu (q,v)  = 0 \} 
\end{eqnarray}
where
\begin{definition} The   horizontal space at $q$  is the  orthogonal complement (rel. the mass metric) of the sum of first two subspaces $\R q$ and $i \R q$,
i.e it is the orthogonal  complement to the  $\C$-span of $q$. 
\end{definition}
\noindent  Refer to exercise \ref{group_perp} and the definition of $\nu$ to see why the horizontal space at $q$ is, as described above,
the zero locus of $J(q,v)$ and $\nu(q, v)$.  

Unit vectors spanning the scale and rotation    spaces  are   $s$ and $is$.
Consequently, if  we take a $v \in T_q \E_{cm}$ and decompose it accordingly we get
\begin{equation}
\label{Saari}
v = \langle s, v \rangle s + \langle is, v \rangle is + v_{hor}
\end{equation}
and the scale invariant version:
\begin{equation}
\label{Saari2}
y  = \nu  s + \tilde J  is + y_{hor} ; \qquad \nu = \langle s, y \rangle, \tilde J = \langle is, y \rangle
\end{equation}
where the subscript ``hor'' on $v$ and $y$ denote their   orthogonal projections onto the horizontal subspace.

{\sc Remark}  D. Saari pointed out the importance to celestial mechanics  of the   horizontal-vertical splitting of eq (\ref{Saari} )
 and hence this splitting is often called the ``Saari decomposition''.

\subsubsection{Proof of proposition \ref{blowupQuotient}}
The decomposition (eq (\ref{Saari2})) of $y$  is $S^1$-equivariant.  The coefficients of the first two terms
$\nu$ and  $\tilde J = \langle is ,y \rangle$ are $S^1$-invariant functions and so
are well defined functions on the quotient ${\mathcal P}_N$.  
The horizontal term $y_{hor}$ , as $y$ varies at fixed $s$, sweeps out the   horizontal subspace at $s$ and these subspaces, as $s$ varies, 
forms the horizontal distribution associated to a   connection on the principal $S^1$-bundle $S_{cm} \to \C \bP^{N-2}$. 
  It is a basic fact
about principal $G$-bundles with connection that   the union of the  horizontal spaces for the connection  forms a
$G$-equivariant  vector bundle
over the total space, and the  quotient of  this  vector bundle by $G$  
is canonically isomorphic to the tangent space to the base space. Writing $[s,y]$ to denote the $S^1$-equivalence class of
the pair $(s,y)$ we see that  the set of all $[s, y_{hor}]$'s  forms $T \C \bP^{N-2}$. Now $s$, together with $(y_{hor},  \nu, \tilde J)$  determine $y$
uniquely.  It follows  that the map $[s,y] \mapsto ([s, y_{hor}], (\nu, \tilde J))$ is a vector bundle    isomorphism between 
the vector bundles   $(S_{cm} \times \E_{cm})/S^1$ and $T \C \bP ^{N-2} \times \R^2$ over $\C \bP ^{N-2}$.  
The radial scaling coordinate  $r$   ``goes along for the ride'' without change.
QED
\vskip .4cm

Because  the decompositions of equations (\ref{Saari}, \ref{Saari2}) are orthogonal and the second decomposition
is scale invariant it follows that total    kinetic energy   decomposes as
\begin{equation}
\label{KEdecomp}
\begin{aligned}
K(q,v)  &= &\frac{1}{2} \frac{\nu ^2}{r}  +  \frac{1}{2} \frac{J^2}{r^2} +  \frac{K_{shape} ([s,y_{hor}])}{r}\\
            &=& \frac{1}{r} (\frac{\nu ^2}{2}  +   \frac{\tilde J^2}{2} +  K_{shape} )
\end{aligned}
\end{equation}
The final  term $K_{shape}$ is formed by computing the squared length of the horizontal factor $y_{hor}$
and is canonically identified with the   kinetic energy of   the standard (Fubini-Study) metric on  
the shape space $\C \bP ^{N-2}$.  

{\sc Remark.} The kinetic energy decomposition (\ref{KEdecomp}) shows that
for  $J \ne 0$ the manifolds $M^{int} (H_0, J)$ is already  closed in ${\mathcal P}_N$
so that 
\begin{equation}M(h, J) = M^{int}(h, J)
\label{Jne0closed}
\end{equation}
  Indeed, the energy equation $rh = \tilde H$  shows that $\tilde U \ge \frac{1}{2} J^2 / r + O(r)$ holds on 
  $M^{int} (h, J)$ which shows that if for a sequence   $p_i \in M^{int} (H_0, J)$
  we have  that $r(p_i) \to 0$ then $U(s_i) \to \infty$ so that the shape $s_i$ of these  points $p_i$ are converging to the collision
  locus ${\mathcal C} \subset \C \bP ^{N-2}$ on the shape space.  But we deleted ${\mathcal C}$ in forming ${\mathcal P}_N$.

\subsection{Euler-Lagrange family in reduced coordinates}

We   follow Moeckel and look into  what a planar     central configuration family of section \ref{ELfamily}
such as the Euler or Lagrange family looks like in the coordinates of   ${\mathcal P}_N$.  

Let $s_{cc}$ be a planar  central configuration  and $[s_{cc}] \in \C \bP^{N-2}$  the corresponding point in shape space.
During the evolution of  the associated family, this shape does not change. Only the size $r$ and angle $\theta$
of of the configuration changes. This size and angle change is specificied by $\lambda = \lambda(t)  = r e^{i \theta}$ where $\lambda(t)$
solves the Kepler problem as per exercise \ref{planarccs}.
Since the shape does not change, the  shape velocity $y_{hor}$  is identically zero along each of these solutions and so   $K_{sh} =0$.
Thus along such a solution 
$$\tilde K = \frac{1}{2} \nu ^2 + \frac{1}{2} \tilde J^2 = \frac{1}{2} \nu ^2 + \frac{1}{2} \frac{J^2}{r}$$
(see eq. (\ref{KEdecomp}))
and the only variables which change  are $(r, \nu, \tilde J)$ among the full set of variables  $(r, [s, y_{hor}], (\nu, \tilde J))$ of  ${\mathcal P}_N =  [0,\infty) \times T (\C \bP^{N-2} \setminus {\mathcal C} ) \times \R^2$
But $\tilde J = r^{1/2} J$ and $J$ is constant along solutions so the change of $r$ and choice of $J$ determines the change of $\tilde J$.
 So we can think of the only variables being $\nu, r$.  

Fix the energy $h$.  We can then  view the central configuration  family as a one-parameter family of curves in the $(\nu, r)$ plane, the parameter  being  the
angular momentum $J$. Indeed the energy equation reads: 
$$r h = \frac{1}{2} \nu ^2 + \frac{1}{2} J^2/ r  - U (s_{cc}).$$
and since $U(s_{cc})$ is constant, this defines a one-parameter family of curves. 
We plot these curves in the $\nu, r$ plane for various values of the angular momentum $J$ below in figure \ref{fig:ccplots}.
\begin{figure}[h]
\label{fig:ccplots}
\center
\scalebox{0.5}{\includegraphics{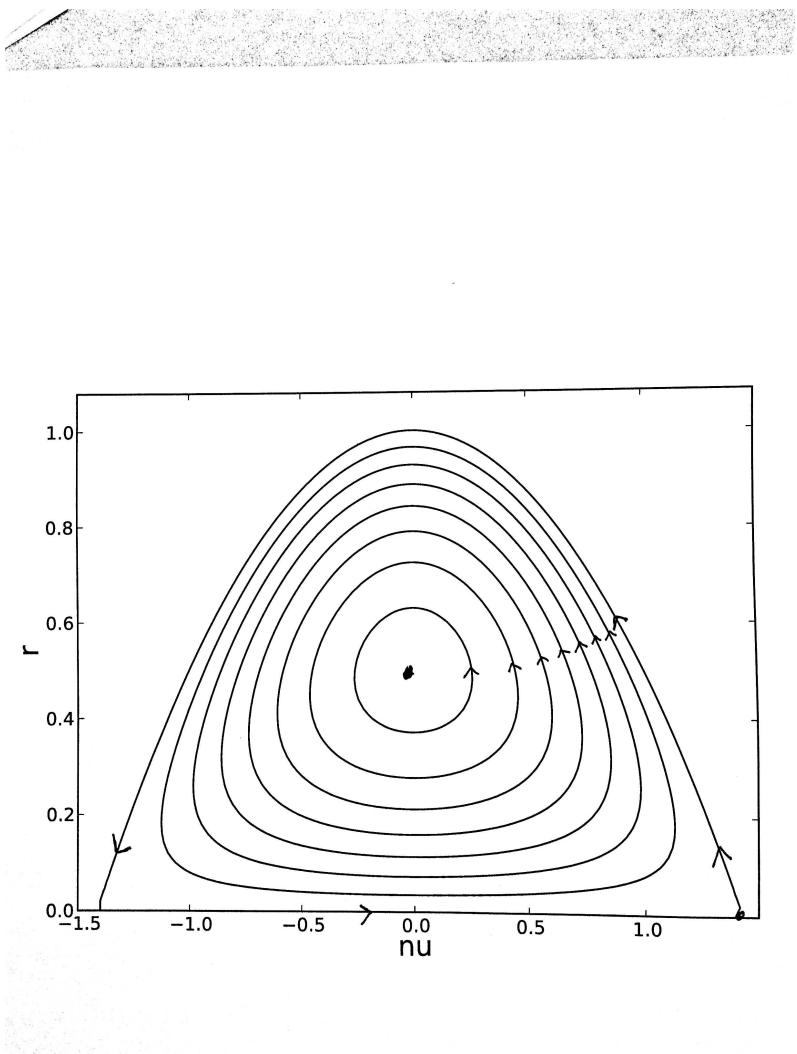}}
\caption{A central configuration family in $\nu, r$ coordinates. The arch and `floor' $r =0$ comprise the rest cycle}  
\end{figure}

Observe the rest point cycle in this picture:  the closed curve passing thru
the two equilibria.  This curve is the union of two solution curves, a top arch which is an honest solution, and a bottom return curve.
The top arch is the
ejection-collision orbit first described when we described Lagrange's solution: it explode out of 
total collision along the shape $s_{cc}$ achieves a   maximum size and shrink back to triple collision. It connects
the rest point $s \in C$ having shape $s_{cc}$ and $\nu = \sqrt{2 U(s_{cc})} > 0$
the rest point $s_*$ having shape $s_{cc}$ and $\nu = - \sqrt{2 U(s_{cc})} < 0$.
This top arch lies on $M(h, 0)$.  The bottom `return road' lies on the full  collision manifold $\{r=0, \tilde H =0\}$
and yields a return route from $s_*$ to $s$.  This  rest point cycle is the  
 limit of the family of the  periodic central configuration solutions  with $J \ne 0$  as $J \to 0$. 
 
 {\bf Notational Convenience.}  
 We have just  used  the   symbol $M(h,0) \subset {\mathcal P}_N$   for what used to be
 a submanifold of the phase space {\it before quotient}.    We will continue to use the same
 notation for any $G$ invariant  submanifold or function on phase space
before or after the quotient procedure.  Thus we have:
 $$C, M_0, M(h),  M^{int} (h),  M(h,J_0),  \text{ etc. }   \subset  {\mathcal P}_N.$$

\section{A gradient (like) flow!}  

The   dominant   aspect of the flow on the full collision manifold $M_0$ is that $-\nu$ acts like a Liapanov function.
\begin{exer}  Use equations (\ref{blowup}) to  derive the identity
\begin{equation}  \nu ^{\prime} = \tilde K - \frac{1}{2} \nu^2 + \tilde H
\label{nuprimeexer}
\end{equation}
(See for example Moeckel \cite{M_near},  eq. (1.6).)
\end{exer}
\begin{exer}  Use the ``Saari decomposition'' of kinetic energy (eq (\ref{KEdecomp})) to show that
$$\tilde K-  \frac{1}{2} \nu^2   = K_{sh} + \frac{1}{2} \tilde J^2.$$
Conclude, using the previous exercise,  that
$$\nu^{\prime} = K_{sh}+ \frac{1}{2} \tilde J^2  \ge 0  
\text{ on } M_0 = \{ r =0, \tilde H = 0 \}.$$
\end{exer}

You have proved much of
\begin{proposition} 
\label{gradlike}$\nu^{\prime} \ge 0$ everywhere on the full collision manifold $M_0$.
Moreover $\nu$ is constant along a solution lying in  $M_0$ if and only if that solution
is one of the equilibria. 
\end{proposition} 

{\bf Remark.} A flow  is called  ``gradient-like'' if it admits a continuous  function $f$ which is strictly monotone decreasing
along all   solution curves except   equilibria.    (See   Robinson \cite{Robinson} p. 357.)
The proposition thus asserts that  the blown-up flow is gradient-like on the full collision manifold $M_0$  relative to the function $f= -\nu$.
\vskip .4cm

{\sc Proof of proposition \ref{gradlike}.}
In the exercise you proved that   $\nu^{\prime} $ is positive
everywhere except at  the  points where it is zero.
We must then  show that any solution which lies on the locus $\nu^{\prime} =0$
is an equilibrium.  We see that
$\nu^{\prime} = 0$ if and only if $K_{sh} = 0 = \tilde J$.  
Now  
 $d \tilde J / d \tau = -\frac{1}{2} \nu \tilde J$.
 (This  holds both on and off the collision manifold).
 It follows that any solution starting on   $\tilde J = 0$ remains on the locus   $\tilde J = 0$. 
   $K_{sh} ([s,y]) = 0$ if and only if   $y_{hor} = 0$ in which case,
 both the $y_{hor}$ and the $is$ term (from $\tilde J = \langle is , y \rangle$) in the decomposition of $y$
 are zero so that   
 $y = \lambda s $ with $\lambda \in \R$.
Take inner products with $s$ to get   $\lambda = \nu$. 
 Now assume we have a solution curve $(s(\tau), y(\tau)$ lying on
 the locus $\tilde J = 0, K_{sh} = 0$.  Differentiating  the equation $y(\tau) = \nu(\tau) s(\tau)$
 using the blow-up equations we see that  $y^{\prime} = \nu^{\prime} s + \nu s^{\prime}$.
 But $\nu^{\prime} =0$ by assumption and $s^{\prime} = y - \nu s= 0$ by the blow-up equations, 
 so $y^{\prime} =0$ along the solution:  our curve is an equilibrium.
 
 QED
 
Pause for a moment. 
 Reflect how different flow on the full collision locus is from a Hamiltonian flow
 on an energy level set.

\subsection{Moeckel's manifold with corner into a manifold with a T}

\begin{figure}[h]
\label{T}
\center
\scalebox{0.5}{\includegraphics{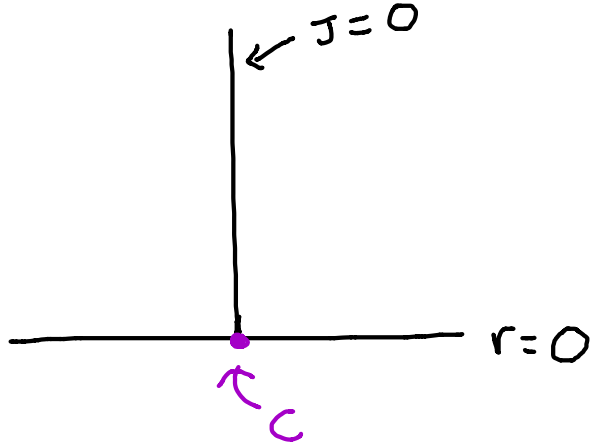}}
\caption{$\hat M (h)$ inside $M(h)$ is the zero level set of $r \tilde J$.}  
\end{figure}

In \cite{Moeckel1}, at the beginning of section 2,  Moeckel constructs  a certain manifold  
with corners  in preparation
for perturbing the heteroclinic tangles lying on $M(h, 0)$ into the realms of
$M(h, \epsilon)$.   (He denotes his  manifold with a corner by  $M_{0+}$ and later simply $M$.)  
Dynamics on this manifold-with-corners is essential to our proof of theorem 1.   
 I had a hard time making sense of this manifold.  
 I     rederive  what Moeckel  did in a   slightly different way.  I get 
 a  ``manifold with a T'' instead of Moeckels manifold with a corner.
  A  ``T '' is made out of two corners, or ``L''s (one reflected relative to the other)    joined along their vertical edge.  One of these  corners  is Moeckel's
manifold with a corner and the other is a reflection of it.  The vertex, or corner itself,  is  our good friend $C$ the standard collision manifold.
  (Figure \ref{T} .)
 
Recall that $M(h)$ is a hypersurface in ${\mathcal P}_N$, and as such is a  manifold
with boundary, whose boundary is  our friend full collision manifold $M_0 = \{r = 0, \tilde H = 0\}$. 
\begin{definition}    $\hat M (h) = M(h,0 ) \cup  M_0  \subset \hat M(h)$.
\end{definition}
$\hat M (h)$  is a codimension 1   subvariety of the smooth manifold with boundary $M(h)$.  
It is the zero locus of the function $r \tilde J$ restricted to $M(h)$ and as such has   two algebraic components
:  $r= 0$ which is our full collision manifold $M_0$, and $J = 0$ which forms $M(h,0)$. The singular locus of $\hat M(h)$
is the intersection $C = \{r = 0 , \tilde J =0 \}$ of these two components.   All the  rest point cycles described above  associated to the central configurations
lie on this $\hat M( h)$.   $\hat M (h)$ is comprised of two ``manifolds with corners'', namely $\{r \tilde J =0 , \tilde J \ge 0 \}$ and
 $\{r \tilde J =0 , \tilde J \le 0 \}$.  The first of these is Moeckel's.

$\hat M(h)$ is to be viewed as the limit as $J\to 0$
of the manifolds $M(h, J)$. 

\begin{proposition}  For $S  \subset \R$ a subset of the line of angular momentum values, 
set $M^{int} (h, S) = \cup_{J \in S} M^{int} (h, J)$.  Then
$\hat M(h) = \cap_{\epsilon > 0 }  M^{int} (h, (-\epsilon, \epsilon))$.
\end{proposition}

The proof of the proposition follows in a routine way from our expressions for energy, from $rh = \tilde H$, $J = r^{-1/2} \tilde J$ and 
the kinetic energy decomposition of eq \ref{KEdecomp}.
It is useful to recall, eq (\ref{Jne0closed}) that the  $M(h, J) = M^{int}(h, J)$ are closed for $J \ne 0$. 

As an alternative to the description of the proposition,    we can either let $J \to 0$ from   above or below.
Set
$$\hat M_+ (h) = \lim_{J \to 0^+}  M(h, J)$$
and 
$$\hat M_- (h) =\lim_{J \to 0^-}  M(h, J) .$$  
Then one can show without difficulty that 
$$\hat M (h) = M_+ (h) \cup M_- (h), $$
with $\hat M_{+} (h) = \{ p \in \hat M(h), \tilde J \ge 0 \}$
and $\hat M_{-} (h) = \{ p \in \hat M(h), \tilde J \le 0 \}$ being the 
two manifolds with corners described earlier,  Moeckel's manifold with a corner being $\hat M_+$.

{\it What is a manifold  with a `T'?}   
Suppose we have two real-valued  functions $x, y$ on an $n$-dimensional manifold  $Q$ such that    $0$ is a regular value 
for both functions and  $(0,0)$ is a regular value
of the map $(x,y): Q \to \R^2$.  Then the locus $\{ xy = 0,  y \ge 0 \}$
is a manifold with a T.    Its singular locus is $\{x = y = 0 \}$. 
A manifold with a T is  locally  diffeomorphic to the product of the ``upside down T''
  $xy = 0, y \ge 0$ in the xy plane,  by an $\R^{n-2}$.  
  See figure \ref{T}.

\subsection{Finishing up the proof of theorem 1.} 
The  idea of Moeckel is that hyperbolic structures   persist on perturbation, and 
that the  various stable-unstable connections
between Euler and Lagrange central configuration points  on on $\hat M (h)$   are sufficiently ``hyperbolic'' that they   persist into $M(h, \epsilon)$
for $\epsilon \ne 0$   small.  Nonzero angular momentum is needed to   get orbits connecting from $R$'s to $R^*$
in finite time since the rest cycle of figure \ref{fig:ccplots} takes infinite blown-up time.   Moeckel cannot carry out the `perturbation of hyerbolic''  idea literally  because he
cannot establish the needed  hyperbolicity or transversality.  Instead, following  an earlier idea of Easton, he
replaces   hyperbolicity by a weaker  notion of 
``topologically transverse''  between collections of ``windows'' transverse to the flow.
This notion is sufficiently flexible and stable to allow Moeckel to perturb the various formal
connections to get actual orbits realizing walks in the abstract graph introduced in section \ref{3}. 
By following the details of his proof, three decades later,  we were able to verify that his realizing solutions
when projected onto the shape space do indeed stay $C^0$-close to the concrete connection graph as described in section \ref{3}.

The hypothesis of equal or near equal masses is needed to insure that (some of) the eigenvalues for the linearization at the Euler equilibria
are complex.  This complexity implies a ``spiralling'' of the Lagrange stable/unstable manifolds around the Euler unstable/stable manifolds and is needed
to insure that all connections in the abstract connection graph are realized. 

\section{Acknowledgements}

I am grateful to acknowledge critique and feedback in the preparation of this article from Alain Albouy, Alain Chenciner,
and Rick Moeckel.  Wu-yi Hsiang asked the central question which inspired this research. 
A conversation with Carl\'es Simo   completely 
redirected my methods to the ones that were ultimately successful.  The whole project owes its existence to Rick Moeckel. I am thankful to the participants of the CIMAT school in Guanajuato, Mexico, 
 to the organizer of that school Rafael Herrera, to Gil Bor for good Israeli salads, company,  and piano playing and
 to Patricia Carral and Eyal Bor for  their hospitality and use of a comfortable  bed during this school.  Finally I wish to thank NSF grant DMS-1305844 for essential support.

\


\

\end{document}